\documentclass[12pt]{article}
\usepackage{color}
\definecolor{darkblue}{rgb}{0.00,0.25,0.50}
\usepackage[colorlinks,filecolor=blue,citecolor=darkblue]{hyperref}

\usepackage{amsfonts,amssymb,amsmath,amsthm}
\usepackage{url}
\usepackage{enumerate}
\usepackage[ukrainian, russian, english]{babel}
\usepackage[cp1251]{inputenc}
\begin{document}\selectlanguage{ukrainian}
\thispagestyle{empty}

\title{}

\begin{center}
\textbf{\Large Порядкові оцінки найкращих наближень та наближень сумами Фур'є в рівномірній метриці
 класів згорток періодичних функцій невеликої  гладкості }
\end{center}
\vskip0.5cm
\begin{center}
А.~С.~Сердюк${}^1$, Т.~А.~Степанюк${}^2$\\ \emph{\small
${}^1$Інститут математики НАН
України, Київ\\
${}^2$Східноєвропейський національний університет імені Лесі
Українки, Луцьк\\}
\end{center}
\vskip0.5cm


\begin{abstract}
Получены точные по порядку оценки наилучших равномерных приближений и равномерных приближений суммами Фурье классов сверток периодических функций, принадлежащих единичным шарам пространств
 $L_{p}, \ {1< p<\infty}$, с производящим ядром $\Psi_{\beta}$, абсолютные величины $\psi(k)$ коэффициентов Фурье которого таковы, что  $\sum\limits_{k=1}^{\infty}\psi^{p'}(k)k^{p'-2}<\infty$,  $\frac{1}{p}+\frac{1}{p'}=1$, а произведение $\psi(n)n^{\frac{1}{p}}$ не может стремится к нулю быстрее степенных функций.

\vskip 0.5cm

We obtain exact for order estimates of best uniform approximations and  uniform approximations
by Fourier sums of classes of convolutions the periodic functions   belong to unit balls of spaces  $L_{p}, \ {1\leq p<\infty}$,
with generating kernel $\Psi_{\beta}$, whose absolute values of Fourier coefficients $\psi(k)$ are such that  $\sum\limits_{k=1}^{\infty}\psi^{p'}(k)k^{p'-2}<\infty$,  ${\frac{1}{p}+\frac{1}{p'}=1}$, and product $\psi(n)n^{\frac{1}{p}}$ can't tend to nought faster than power functions.
\end{abstract}


Позначимо через $C$ --- простір $2\pi$--періодичних неперервних функцій, у
якому норма задана за допомогою рівності
$$
{\|f\|_{C}:=\max\limits_{t}|f(t)|};$$
$L_{p}$, $1\leq p\leq\infty$, --- простір $2\pi$-періодичних
сумовних функцій $f$ зі скінченною нормою $\| f\|_{p}$, де
$$
\|f\|_{p}:={\left\{\begin{array}{cc}
\Big(\int\limits_{-\pi}^{\pi}|f(t)|^{p}dt\Big)^{\frac{1}{p}}, \ & 1\leq p<\infty, \\
\mathop{\rm{ess}\sup}\limits_{t}|f(t)|, &
p=\infty. \
  \end{array} \right.}
$$

Нехай $f$ --- функція із $L_{1}$, ряд Фур'є  якої має вигляд
$$
\frac{a_{0}}{2}+\sum_{k=1}^{\infty}(a_{k}\cos kx+
b_{k}\sin kx).
$$
Нехай, далі, $\psi(k)$ --- довільна фіксована послідовність дійсних чисел і $\beta$
--- фіксоване  дійсне число. Тоді якщо ряд
$$
\sum_{k=1}^{\infty}\frac{1}{\psi(k)}\Big(a_{k}\cos\Big(kx+\frac{\beta\pi}{2}\Big)
+b_{k}\sin\Big(kx+\frac{\beta\pi}{2}\Big)\Big)
$$
\noindent є рядом Фур'є деякої сумовної функції $\varphi$, то цю
функцію  називають (див., наприклад,  \cite[с.
132]{Stepanets1}) $(\psi,\beta)$-похідною функції $f$ і позначають через
$f_{\beta}^{\psi}$.
Множину функцій $f$, у яких існує $(\psi,\beta)$-похідна
позначають через $L_{\beta}^{\psi}$.

Покладемо
$$
B_{p}^{0}:=\left\{\varphi\in L_{p}: \ ||\varphi||_{p}\leq 1,
\ \varphi\perp1\right\}, \ \ 1\leq p\leq\infty.
$$
Якщо $f\in L^{\psi}_{\beta}$, і водночас $f^{\psi}_{\beta}\in
B_{p}^{0}$, то кажуть, що функція $f$ належить
класу
  $L^{\psi}_{\beta,p}$.
  Позначимо також
  $
    C^{\psi}_{\beta}=C\cap L^{\psi}_{\beta}, \ \ C^{\psi}_{\beta,p}=C\cap L^{\psi}_{\beta,p}.
  $

 Будемо розглядати послідовності $\psi(k)$ такі, що $\psi(k)k^{\frac{1}{p}}$ монотонно незростає і  $\sum\limits_{k=1}^{\infty}\psi^{p'}(k)k^{p'-2}<\infty$, $1<p<\infty$, ${\frac{1}{p}+\frac{1}{p'}=1}$. Тоді з урахуванням  леми 12.6.6 монографії  \cite[с. 193]{Zigmund2}
  та твердження 3.8.3 монографії \cite[с.
139]{Stepanets1} функції $f$ з множини
$C^{\psi}_{\beta,p}, \ 1\leq p\leq \infty$,
 для всіх
$x\in\mathbb{R}$ зображуються   за допомогою згортки
\begin{equation}\label{conv}
f(x)=\frac{a_{0}}{2}+\frac{1}{\pi}\int\limits_{-\pi}^{\pi}\Psi_{\beta}(x-t)\varphi(t)dt,
  \ \varphi\in B^{0}_{p},
\end{equation}
де
\begin{equation}\label{kernel}
\Psi_{\beta}(t)=\sum\limits_{k=1}^{\infty}\psi(k)\cos
\big(kt-\frac{\beta\pi}{2}\big), \ \  \beta\in
    \mathbb{R},
\end{equation}
і при цьому майже скрізь $\varphi=f^{\psi}_{\beta}$.

 При $\psi(k)=k^{-r}$, $r>0$, ядра $\Psi_{\beta}(t)$ вигляду (\ref{kernel}) є ядрами Вейля--Надя
 $B_{r,\beta}(t)$
 \begin{equation}\label{brb}
B_{r,\beta}(t)=\sum\limits_{k=1}^{\infty}k^{-r}\cos
\big(kt-\frac{\beta\pi}{2}\big), \ \  \beta\in
    \mathbb{R},
\end{equation}
а класи  функцій $f$, що зображуються у вигляді
\begin{equation}\label{cc}
f(x)=\frac{a_{0}}{2}+\frac{1}{\pi}\int\limits_{-\pi}^{\pi}B_{r,\beta}(x-t)\varphi(t)dt,
  \ \varphi\in B^{0}_{p},
\end{equation}
 є  відомими класами
Вейля-Надя  $W^{r}_{\beta,p}, \ 1<p<\infty$.
Зрозуміло, що при $r>\frac{1}{p}$, ${1<p\leq\infty}$, для довільних $\beta\in  \mathbb{R}$ має місце включення $W^{r}_{\beta}\subset C$.

Позначимо через ${\mathfrak M}$ множину усіх опуклих донизу, неперервних функцій $\psi(t)$, $t\geq 1,$ таких, що $\lim\limits_{t\rightarrow\infty}\psi(t)=0$. І будемо вважати, що послідовність $\psi (k),\ k\in \mathbb{N}$, яка задає клас $C^{\psi}_{\beta,p}, \ 1\leq p\leq \infty$,  є звуженням на множину
натуральних чисел
функцій $\psi(t)$ із ${\mathfrak M}$.

Вслід за О.І. Степанцем (див., наприклад,  \cite[с.
160]{Stepanets1}),
 за допомогою характеристики $\mu(\psi;t)$ функцій $\psi$ із $\in{\mathfrak M}$ вигляду
\begin{equation}\label{mu}
\mu(t)=\mu(\psi;t):=\frac{t}{\eta(t)-t},
\end{equation}
де $\eta(t)=\eta(\psi;t):=\psi^{-1}\left(\psi(t)/2\right)$, $\psi^{-1}$ --- обернена до $\psi$ функція, з множини ${\mathfrak M}$ виділимо наступні підмножини:
\begin{equation}\label{m0}
\mathfrak{M}_{0}=\left\{\psi\in \mathfrak{M}: \ \ \exists K>0 \ \ \ \  \forall t\geq1 \ \ \ \
0< \mu(\psi;t)\leq K<\infty \right\},
\end{equation}
\begin{equation}\label{mc}
\mathfrak{M}_{C}=\left\{\psi\in \mathfrak{M}: \ \ \exists K_{1}, K_{2}>0 \ \ \ \ \forall t\geq1 \ \ \ \
K_{1}\leq \mu(\psi;t)\leq K_{2}<\infty \right\},
\end{equation}
$$
\mathfrak{M}^{+}_{\infty}=\left\{\psi\in \mathfrak{M}: \ \ \
\mu(\psi;t)\uparrow\infty \right\}.
$$
 В (\ref{m0}) i (\ref{mc}) $K, \ K_{1}, \ K_{2}$, взагалі кажучи, можуть залежати від $\psi$.
Очевидно, що ${\mathfrak{M}_{C}\subset\mathfrak{M}_{0}}$.
Відмітимо, що природними представниками множини $\mathfrak{M}_{C}$ є функції $\psi(t)=t^{-r}$, $r>0$,
 множини $\mathfrak{M}_{0}\setminus \mathfrak{M}_{C}$ --- функції $\psi(t)=\ln^{-\varepsilon}(t+1), \ \varepsilon>0$, а  множини
 $\mathfrak{M}^{+}_{\infty}$ --- функції $e^{-\alpha t^{r}}$, $r>0, \ \alpha>0$.

Крім величини (\ref{mu}) для функцій $\psi\in\mathfrak{M}$ важливу роль відіграє  характеристика
\begin{equation}\label{for301}
\alpha(\psi;t):=\frac{\psi(t)}{t|\psi'(t)|}, \ \ \psi'(t):=\psi'(t+0).
\end{equation}

В  \cite[с.
160]{Stepanets1} було доведено, що необхідною і достатньою умовою належності
 функції $\psi\in\mathfrak{M}$ до множини $\mathfrak{M}_{0}$ є умова
$$
\alpha(\psi;t)>K>0 \ \forall t\geq1;
$$
а необхідною і достатньою умовою того, щоб функція $\psi\in\mathfrak{M}$ належала до множини $\mathfrak{M}_{C}$ є умова
$$
0<K_{1}\leq\alpha(\psi;t)\leq K_{2}<\infty \  \ \ \ \forall t\geq1.
$$

Якщо $\psi\in\mathfrak{M}^{+}_{\infty}$, то (див., наприклад,  \cite[с.
97]{Step monog 1987}) множини
$C^{\psi}_{\beta}$ складаються з нескінченно диференційовних
функцій. З іншого боку, як показано в  \cite[с.
1692]{Stepanets_Serdyuk_Shydlich}, для кожної нескінченно
диференційовної $2\pi$--періодичної функції $f$ можна вказати функцію $\psi$ з множини
$\mathfrak{M}^{+}_{\infty}$ таку, що $f\in C^{\psi}_{\beta}$ для
довільних $\beta\in\mathbb{R}$.

Для класів  $C^{\psi}_{\beta,p}$  будемо розглядати величини
$$
{\cal E}_{n}(C^{\psi}_{\beta,p})_{C}=\sup\limits_{f\in
C^{\psi}_{\beta,p}}\|f(\cdot)-S_{n-1}(f;\cdot)\|_{C}, \ \ 1\leq p\leq\infty,
$$
де $S_{n-1}(f;\cdot)$ --- частинні суми Фур'є порядку $n-1$ функції $f$, а також
 найкращі рівномірні наближення класів
$C^{\psi}_{\beta,p}$ тригонометричними
поліномами  порядку не вищого за ${n-1}$, тобто величини вигляду
$$
{E}_{n}(C^{\psi}_{\beta,p})_{C}=\sup\limits_{f\in
C^{\psi}_{\beta,p}}\inf\limits_{t_{n-1}\in\mathcal{T}_{2n-1}}\|f(\cdot)-t_{n-1}(\cdot)\|_{C},
\ 1\leq p\leq \infty,
$$
де $\mathcal{T}_{2n-1}$ --- підпростір усіх тригонометричних
поліномів $t_{n-1}$ порядку не вищого за ${n-1}$.

В роботі розв'язується задача про  знаходження точних порядкових оцінок для
 величин ${\cal E}_{n}(C^{\psi}_{\beta,p})_{C}$ i ${ E}_{n}(C^{\psi}_{\beta,p})_{C}$ при $1< p<\infty$, $\beta\in\mathbb{R}$.

Для класів Вейля--Надя   $W^{r}_{\beta,p}$, $r>\frac{1}{p}$,
$\beta\in\mathbb{R}$, ${1\leq p\leq\infty}$, точні порядкові оцінки
величин ${\cal E}_{n}(W^{r}_{\beta,p})_{C}$ i
${E}_{n}(W^{r}_{\beta,p})_{C}$ відомі (див., наприклад,  \cite[с.
47--49]{T}.
Крім того,  для величин ${\cal
E}_{n}(W^{r}_{\beta,\infty})_{C}$, $r>0$, $\beta\in \mathbb{R}$ при $n\rightarrow\infty$ відомі асимптотичні рівності
(див., наприклад,  роботи  \cite{Kol}--\cite{Nik2}), а для величин найкращих наближень ${E}_{n}(W^{r}_{\beta,\infty})_{C}$  при усіх
$n\in\mathbb{N}$, встановлено їх точні значення
(див.  роботи
\cite{Nik2}--\cite{SU}).

 Для класів  $C^{\psi}_{\beta,p}$ при $p=2$ і $\beta\in \mathbb{R}$ за умови
$\sum\limits_{k=1}^{\infty}\psi^{2}(k)<\infty$ в роботі
\cite{Serdyuk} було доведено рівність
$$
{\cal E}_{n}(C^{\psi}_{\beta,2})_{C}=\frac{1}{\sqrt{\pi}}\Big(\sum\limits_{k=n}^{\infty}\psi^{2}(k)\Big)^{\frac{1}{2}}, \ n\in\mathbb{N}.
$$

В \cite{Serdyuk_grabova}
 у випадку, коли
$\psi\in B\cap\Theta_{p}$, де $\Theta_{p}$, $1\leq p<\infty$, ---
множина незростаючих функцій $\psi(t)$, для яких існує
стала $\alpha>\frac{1}{p}$ така, що функція $t^{\alpha}\psi(t)$
майже спадає (тобто знайдеться додатна стала $K$ така, що
$t^{\alpha}_{1}\psi(t_{1})\leq
Kt^{\alpha}_{2}\psi(t_{2})$ для будь--яких $t_{1}>t_{2}\geq 1$), а $B$ --- множина
незростаючих додатних функцій $\psi(t)$, $t\geq 1$, для кожної з
яких можна вказати додатну сталу $K$ таку, що $
\frac{\psi(t)}{\psi(2t)}\leq K, \ \  t\geq 1 $, показано, що існують додатні величини $K^{(1)}$ і $K^{(2)}$,
які можуть  залежати лише  від $\psi$ і $p$ такі, що для довільних $1<p<\infty$, $\beta\in \mathbb{R}$ i $n\in \mathbb{N}$
 виконуються співвідношення
\begin{equation}\label{gr_ser}
K^{(2)}\psi(n)n^{\frac{1}{p}}\leq{
E}_{n}(C^{\psi}_{\beta,p})_{C}\leq{\cal
E}_{n}(C^{\psi}_{\beta,p})_{C}\leq K^{(1)}
\psi(n)n^{\frac{1}{p}}.
\end{equation}

У випадку $\psi\in \mathfrak{M}_{\infty}^{+}$ порядкові оцінки величин ${\cal E}_{n}(C^{\psi}_{\beta,p})_{C}$
і ${ E}_{n}(C^{\psi}_{\beta,p})_{C}$, $1\leq p<\infty$,  знайдені в роботах \cite{Step monog 1987}, \cite{Rom}, \cite{S_S}.

Мета даної роботи полягає у знаходженні двосторонніх  оцінок для
 величин ${\cal E}_{n}(C^{\psi}_{\beta,p})_{C}$ та ${ E}_{n}(C^{\psi}_{\beta,p})_{C}$  у випадку коли  функція
 $g_{p}(t)=\psi(t)t^{\frac{1}{p}}$ належить до множини $\mathfrak{M}_{0}$ i $\sum\limits_{k=1}^{\infty}\psi^{p'}(k)k^{p'-2}<\infty$, $1<p<\infty$, $\frac{1}{p}+\frac{1}{p'}=1$. При цьому константи в отриманих оцінках будуть виражені через параметри  класів в явному вигляді.

{\bf Теорема 1.} {\it Нехай $\psi(t)t^{\frac{1}{p}}\in\mathfrak{M}_{0}$ i $\sum\limits_{k=1}^{\infty}\psi^{p'}(k)k^{p'-2}<\infty$, $1<p<\infty$, $\frac{1}{p}+\frac{1}{p'}=1$.
  Тоді   для довільних $n\in \mathbb{N}$ i $\beta\in\mathbb{R}$  мають місце співвідношення
   \begin{equation}\label{theorem_1}
 K_{\psi,p}^{(1)}  \Big(\sum\limits_{k=n}^{\infty}\psi^{p'}(k)k^{p'-2}\Big)^{\frac{1}{p'}} \leq{\cal E}_{n}(C^{\psi}_{\beta,p})_{C}
\leq
K_{\psi,p}^{(2)}\Big(\sum\limits_{k=n}^{\infty}\psi^{p'}(k)k^{p'-2}\Big)^{\frac{1}{p'}},
\end{equation}
в яких $K^{(1)}_{\psi,p}$ і $K^{(2)}_{\psi,p}$ --- додатні величини,
 що   залежать лише  від $\psi$ і $p$.  }

{\bf Доведення теореми 1.}
Згідно з інтегральним зображенням (\ref{conv}), за виконання умов теореми, для довільної
функції $f\in C^{\psi}_{\beta,p}$
 в кожній точці $x\in\mathbb{R}$ справедлива  рівність
\begin{equation}\label{for5}
f(x)-S_{n-1}(f;x)=\frac{1}{\pi}\int\limits_{-\pi}^{\pi}\Psi_{\beta,n}(x-t)f^{\psi}_{\beta}(t)dt, \ \ \|f^{\psi}_{\beta}\|_{p}\leq1, \ f^{\psi}_{\beta}\perp1,
\end{equation}
 де
 \begin{equation}\label{for6}
\Psi_{\beta,n}(t)=
\sum\limits_{k=n}^{\infty}\psi(k)\cos\big(kt-\frac{\beta\pi}{2}\big),  \ \ \beta\in\mathbb{R}, \ \ n\in \mathbb{N}.
\end{equation}

З (\ref{for5}) і  твердження 3.8.1  роботи
 \cite[с. 137]{Stepanets1} одержуємо
 \begin{equation}\label{gg10}
{\cal E}_{n}(C^{\psi}_{\beta,p})_{C}\leq
\frac{1}{\pi}\big\|\Psi_{\beta,n}(\cdot)\big\|_{p'}\|
f^{\psi}_{\beta}(\cdot)\|_{p}\leq\frac{1}{\pi}\big\|\Psi_{\beta,n}(\cdot)\big\|_{p'}, \ 1< p <\infty, \ \frac{1}{p}+\frac{1}{p'}=1.
\end{equation}

Для оцінки величини $\big\|\Psi_{\beta,n}(\cdot)\big\|_{p'}$  буде корисна наступна лема.

{\bf Лема 1.} {\it Нехай  $1<s<\infty$
i $\{a_{k}\}_{k=1}^{\infty}$ --- монотонно незростаюча послідовність додатних чисел така, що $\sum\limits_{k=1}^{\infty}a_{k}^{s}k^{s-2}<\infty$.
Тоді для $L_{s}$--норми функції
$$
h_{\gamma,n}(x)=\sum\limits_{k=n}^{\infty}a_{k}\cos(kx+\gamma),\ \gamma\in\mathbb{R}, \ n\in\mathbb{N},
$$
має місце нерівність
\begin{equation}\label{lemma_1}
\|h_{\gamma,n}(x)\|_{s}\leq \xi(s)\Big(\sum\limits_{k=n}^{\infty}a_{k}^{s}k^{s-2}+a_{n}^{s}n^{s-1}\Big)^{\frac{1}{s}},
 \end{equation}
де
\begin{equation}\label{for100}
\xi(s):=\max\Big\{4\Big(\frac{\pi}{s-1}\Big)^{\frac{1}{s}}, \ \ 14(8\pi)^{\frac{1}{s}} s\Big\}.
\end{equation}
}

{\bf Доведення леми 1.}
Для довільних $n\in\mathbb{N}$, $\gamma\in\mathbb{R}$ i $1<s<\infty$ маємо
$$
\|h_{\gamma,n}(x)\|_{s}^{s}=\int\limits_{-\pi}^{\pi}|h_{\gamma,n}(x)|^{s}dx=
\int\limits_{\frac{\pi}{2n}\leq|x|\leq\pi}|h_{\gamma,n}(x)|^{s}dx+\int\limits_{|x|\leq\frac{\pi}{2n}}|h_{\gamma,n}(x)|^{s}dx=
$$
\begin{equation}\label{for19}
=
J_{s,n}^{(1)}+J_{s,n}^{(2)},
\end{equation}
де
$$
J_{s,n}^{(1)}:=\int\limits_{\frac{\pi}{2n}\leq|x|\leq\pi}|h_{\gamma,n}(x)|^{s}dx,
$$
$$
J_{s,n}^{(2)}:=\int\limits_{|x|\leq\frac{\pi}{2n}}|h_{\gamma,n}(x)|^{s}dx.
$$

Оцінимо величину $J_{s,n}^{(1)}$. Згідно з перетворенням Абеля для довільних $M,N\in\mathbb{N}$
$$
\sum\limits_{k=M}^{N}a_{k}\cos(kx+\gamma)=
$$
\begin{equation}\label{for3}
=\sum\limits_{k=M}^{N-1}\Delta a_{k}\sum\limits_{j=0}^{k}\cos(jx+\gamma)-
a_{M}\sum\limits_{j=0}^{M-1}\cos(jx+\gamma)+a_{N}\sum\limits_{j=0}^{N}\cos(jx+\gamma),
\end{equation}
де $\Delta a_{k}:=a_{k}-a_{k+1}$.

Використовуючи формулу
$$
\sum\limits_{j=0}^{k}\cos(jx+\gamma)=\cos\Big(\frac{k}{2}x +\gamma\Big)\sin
\frac{(k+1)x}{2}\cosec \frac{x}{2},
$$
(див., наприклад,  \cite[с. 43]{Gradshteyn})
та очевидну нерівність
$$
\sin \frac{|x|}{2}\geq\frac{|x|}{\pi}, \ \
 0\leq |x|\leq\pi,
$$
отримаємо
\begin{equation}\label{for4}
 \Big|\sum\limits_{j=0}^{k}\cos(jx+\gamma)\Big|\leq \frac{1}{|\sin\frac{x}{2}|}\leq\frac{\pi}{|x|},  \ \ 0<|x|\leq\pi.
\end{equation}
З умов леми 1 випливає, що $a_{k}\rightarrow 0$ при $k\rightarrow\infty$. Тому
з (\ref{for3}), (\ref{for4}) і того, що внаслідок монотонності  послідовності $a_{k}$ $ \ \Delta a_{k}\geq0$, одержуємо
\begin{equation}\label{for10}
 \Big|\sum\limits_{k=M}^{\infty}a_{k}\cos(kx+\gamma)\Big|\leq\frac{\pi}{|x|}\sum\limits_{k=M}^{\infty}\Delta a_{k}+a_{M}\frac{\pi}{|x|}=a_{M}\frac{2\pi}{|x|},  \ \ 0<|x|\leq\pi, \ \ M\in\mathbb{N}.
\end{equation}
Використавши нерівність (\ref{for10}) при $M=n$, запишемо оцінку
 \begin{equation}\label{for9}
J_{s,n}^{(1)}\leq(2\pi a_{n})^{s}\int\limits_{\frac{\pi}{2n}\leq|x|\leq\pi}\frac{1}{|x|^{s}}dx<\frac{4^{s}\pi}{s-1}a_{n}^{s}n^{s-1}.
\end{equation}
Оцінимо тепер величину $J_{s,n}^{(2)}$. Для довільних  $l,n\in\mathbb{N}, l>n$,
 \begin{equation}\label{for11}
|h_{\gamma,n}(x)|\leq\sum\limits_{k=n}^{l-1}a_{k}+\Big|\sum\limits_{k=l}^{\infty}a_{k}\cos(kx+\gamma)\Big|.
\end{equation}
Покладемо
\begin{equation}\label{eq35}
  A_{n,l}:=\sum\limits_{k=n}^{l-1}a_{k}.
\end{equation}

Беручи до уваги формули (\ref{for11})--(\ref{eq35}), та використавши нерівність (\ref{for10}) при $M=l$, отримаємо
 \begin{equation}\label{for13}
|h_{\gamma,n}(x)|\leq A_{n,l}+a_{l}\frac{2\pi}{|x|}, \ \ 0<|x|\leq\pi.
\end{equation}
Для довільних  $l\geq 2n$, $\frac{\pi}{l+1}<|x|\leq\frac{\pi}{l}$,
$\ l,n\in\mathbb{N}$, в силу монотонного незростання послідовності $a_{k}$, можемо записати
$$
\frac{a_{l}\pi}{|x|}\leq a_{l}(l+1)=a_{l}(l-n+n+1)\leq a_{l}(2(l-n)+1)\leq3a_{l}(l-n)\leq
$$
\begin{equation}\label{for14}
\leq3\sum\limits_{k=n}^{l-1}a_{k}=3A_{n,l}.
\end{equation}
Об'єднуючи (\ref{for13}) і (\ref{for14}), запишемо
\begin{equation}\label{for15}
|h_{\gamma,n}(x)|\leq 7A_{n,l}, \ \frac{\pi}{l+1}< |x|\leq\frac{\pi}{l}, \ \l\geq 2n.
\end{equation}
Із (\ref{for15}) випливає
$$
J_{s,n}^{(2)}=\int\limits_{|x|\leq\frac{\pi}{2n}}|h_{\gamma,n}(x)|^{s}dx=\sum\limits_{l=2n}^{\infty}\int\limits_{\frac{\pi}{l+1}<|x|\leq\frac{\pi}{l}}|h_{\gamma,n}(x)|^{s}dx
\leq2\cdot  7^{s}\sum\limits_{l=2n}^{\infty}\int\limits_{\frac{\pi}{l+1}}^{\frac{\pi}{l}}A_{n,l}^{s}dx\leq
$$
\begin{equation}\label{for16}
\leq
2\pi  7^{s}\sum\limits_{l=2n}^{\infty}\frac{A_{n,l}^{s}}{l^2}\leq2\pi  7^{s}\sum\limits_{l=n+1}^{\infty}\frac{A_{n,l}^{s}}{l^2}.
\end{equation}

При кожному фіксованому $n\in\mathbb{N}$ позначимо через $\alpha_{n}(t), \ t\geq0$, функцію, що означається натупним чином:
$$
\alpha_{n}(t):={\left\{\begin{array}{cc}
a_{k}, \ & k\leq t< k+1, \ k\geq n, \\
0, &
0\leq t<n. \
  \end{array} \right.}
$$
При таких позначеннях має місце рівність ${A}_{n,l}=\int\limits_{0}^{l}\alpha_{n}(t)dt$. Тоді
$$
\sum\limits_{l=n+1}^{\infty}\frac{{A}_{n,l}^{s}}{l^2}=\sum\limits_{l=n+1}^{\infty}\frac{\Big(\int\limits_{0}^{l}\alpha_{n}(t)dt\Big)^{s}}{l^2}
=\sum\limits_{l=n+1}^{\infty}\int\limits_{l}^{l+1}\frac{\Big(\int\limits_{0}^{l}\alpha_{n}(t)dt\Big)^{s}}{l^2}dx\leq
$$
$$
\leq
\sum\limits_{l=n+1}^{\infty}\int\limits_{l}^{l+1}\frac{\Big(\int\limits_{0}^{x}\alpha_{n}(t)dt\Big)^{s}}{(x-1)^2}dx=
\int\limits_{n+1}^{\infty}\frac{\Big(\int\limits_{0}^{x}\alpha_{n}(t)dt\Big)^{s}}{(x-1)^2}dx
=
$$
$$
=
\int\limits_{n+1}^{\infty}\frac{\Big(\int\limits_{0}^{x}\alpha_{n}(t)dt\Big)^{s}}{x^2}dx+
\int\limits_{n+1}^{\infty}\frac{\Big(\int\limits_{0}^{x}\alpha_{n}(t)dt\Big)^{s}}{x^2}\frac{2x-1}{(x-1)^{2}}dx
\leq
$$
$$
=
\int\limits_{n+1}^{\infty}\frac{\Big(\int\limits_{0}^{x}\alpha_{n}(t)dt\Big)^{s}}{x^2}dx+
\Big(\frac{2}{n}+\frac{1}{n^{2}}\Big)\int\limits_{n+1}^{\infty}\frac{\Big(\int\limits_{0}^{x}\alpha_{n}(t)dt\Big)^{s}}{x^2}dx
\leq
$$
\begin{equation}\label{for303}
\leq 4\int\limits_{n}^{\infty}\frac{\Big(\int\limits_{0}^{x}\alpha_{n}(t)dt\Big)^{s}}{x^2}dx.
\end{equation}

Для оцінки останнього інтеграла нам буде корисним наступне твердження, встановлене Харді (див., наприклад,  \cite[с.
40]{Zigmund1}). 

{\bf Лема 2.} {\it Нехай $g(x)$  невід'ємна функція, визначена для $x\geq0$  і нехай $r>1$, $\sigma<r-1$ .
 Тоді, якщо $g^{r}(x)x^{\sigma}$ інтегровна на $(0,\infty)$, то функція $\Big(\frac{1}{x}\int\limits_{0}^{x}g(t)dt\Big)^{r}x^{\sigma}$
 також інтегровна на $(0,\infty)$ і при цьому виконується нерівність
 \begin{equation}\label{Hardi}
\int\limits_{0}^{\infty}\Big(\frac{1}{x}\int\limits_{0}^{x}g(t)dt\Big)^{r}x^{\sigma}dx\leq \Big(\frac{r}{r-\sigma-1}\Big)^{r} \int\limits_{0}^{\infty}g^{r}(x)x^{\sigma}dx.
\end{equation}
 }
Застосовуючи нерівність (\ref{Hardi}) при $g(\cdot)=\alpha_{n}(\cdot)$,  $r=s$, $\sigma=s-2$,  одержуємо
$$
\int\limits_{n}^{\infty}\frac{\Big(\int\limits_{0}^{x}\alpha_{n}(t)dt\Big)^{s}}{x^2}dx\leq
s^{s}\int\limits_{n}^{\infty}\alpha_{n}^{s}(x)x^{s-2}dx=
$$
\begin{equation}\label{for17}
=s^{s}\sum\limits_{l=n}^{\infty}\int\limits_{l}^{l+1}a_{l}^{s}x^{s-2}dx\leq (2s)^{s}\sum\limits_{l=n}^{\infty}a_{l}^{s}l^{s-2}.
\end{equation}
Формули (\ref{for16}), (\ref{for303}) і (\ref{for17}) дозволяють записати оцінку
\begin{equation}\label{for18}
J_{s,n}^{(2)}\leq 8\pi(14s)^{s}\sum\limits_{l=n}^{\infty}a_{l}^{s}l^{s-2}.
\end{equation}
Об'єднавши (\ref{for19}), (\ref{for9})  і (\ref{for18}), маємо
$$
\|h_{\gamma,n}(x)\|_{s}^{s}= J_{s,n}^{(1)}+J_{s,n}^{(2)}\leq
$$
\begin{equation}\label{eqq10}
\leq
\max\Big\{\frac{4^{s}\pi}{s-1} , 8\pi(14s)^{s} \Big\} \Big(
\sum\limits_{l=n}^{\infty}a_{l}^{s}l^{s-2}+a_{n}^{s}n^{s-1}\Big), \ n\in\mathbb{N}, \ 1<s<\infty.
\end{equation}
Із (\ref{eqq10}) випливає (\ref{lemma_1}).
Лему 1 доведено.

Застосуємо до функції $\Psi_{\beta,n}(t)$ вигляду (\ref{for6}) лему 1, поклавши в її умовах
 $a_{k}=\psi(k)$, $\gamma=-\frac{\beta\pi}{2}$,  $s=p'$. Отримаємо оцінку
$$
\|\Psi_{\beta,n}(t)\|_{p'}\leq
 $$
 \begin{equation}\label{eq60}
  \leq \xi(p')\Big(\sum\limits_{k=n}^{\infty}\psi^{p'}(k)k^{p'-2}+\psi^{p'}(n)n^{p'-1}\Big)^{\frac{1}{p'}}, \ 1<p<\infty,
\ \ \frac{1}{p}+\frac{1}{p'}=1,
 \end{equation}
 де характеристика  $\xi(p')$ означається рівністю (\ref{for100}).

Для будь--якої функції $\psi\in\mathfrak{M}$ через $\underline{\alpha}_{n}(\psi)$ i $\overline{\alpha}_{n}(\psi)$, $n\in \mathbb{N}$, позначимо величини
\begin{equation}\label{k}
\underline{\alpha}_{n}(\psi):=\inf\limits_{ t\geq n}\alpha(\psi;t),
\end{equation}
\begin{equation}\label{kk}
\overline{\alpha}_{n}(\psi):=\sup\limits_{ t\geq n}\alpha(\psi;t),
\end{equation}
де характеристика $\alpha(\psi;t)$ означається формулою (\ref{for301}).
В прийнятих позначеннях має місце наступне твердження.

{\bf Лема 3.} {\it Нехай $g_{p}(t):=\psi(t)t^{\frac{1}{p}}$, $\sum\limits_{k=1}^{\infty}\psi^{p'}(k)k^{p'-2}<\infty$, $1<p<\infty$, $\frac{1}{p}+\frac{1}{p'}=1$, $n\in \mathbb{N}$.
  Тоді, якщо $g_{p}\in\mathfrak{M}_{0}$, то виконується нерівність
 \begin{equation}\label{lemma_3}
\psi^{p'}(n)n^{p'-1}\leq\frac{p'}{\underline{\alpha}_{n}(g_{p})}\sum\limits_{k=n}^{\infty}\psi^{p'}(k)k^{p'-2},
\end{equation}

якщо ж  $g_{p}\in\mathfrak{M}_{C}$, то має місце співвідношення}
\begin{equation}\label{lemma_3_1}
\frac{p'}{\overline{\alpha}_{n}(g_{p})}\cdot\frac{n\underline{\alpha}_{n}(g_{p})}{p'+n\underline{\alpha}_{n}(g_{p})}\sum\limits_{k=n}^{\infty}\psi^{p'}(k)k^{p'-2}
\leq\psi^{p'}(n)n^{p'-1}\leq\frac{p'}{\underline{\alpha}_{n}(g_{p})}\sum\limits_{k=n}^{\infty}\psi^{p'}(k)k^{p'-2}.
\end{equation}

{\bf Доведення леми 3.}
Нехай $g_{p}\in\mathfrak{M}_{0}$. Покажемо справедливість нерівності (\ref{lemma_3}).

Оскільки $g_{p}(t)$ монотонно спадає до нуля, то
\begin{equation}\label{eqq2}
\sum\limits_{k=n}^{\infty}\psi^{p'}(k)k^{p'-2}= \sum\limits_{k=n}^{\infty}\frac{g_{p}^{p'}(k)}{k}\geq\int\limits_{n}^{\infty}\frac{g_{p}^{p'}(t)}{t}dt
=\int\limits_{n}^{\infty}\psi^{p'}(t)t^{p'-2}dt.
\end{equation}

Проінтегрувавши частинами останній інтеграл з   (\ref{eqq2}), отримаємо
$$
\int\limits_{n}^{\infty}\psi^{p'}(t)t^{p'-2}dt=
\frac{-1}{p'-1}\psi^{p'}(n)n^{p'-1}-\frac{p'}{p'-1}\int\limits_{n}^{\infty}\psi^{p'-1}(t)\psi'(t)t^{p'-1}dt=
$$
\begin{equation}\label{eqq1}
 =
\frac{-1}{p'-1}\psi^{p'}(n)n^{p'-1}+\frac{p'}{p'-1}\int\limits_{n}^{\infty}\frac{\psi^{p'}(t)t^{p'-2}}{\alpha(\psi;t)}dt.
\end{equation}
З  (\ref{eqq1}) випливає
\begin{equation}\label{eqq3}
\psi^{p'}(n)n^{p'-1} =
p'\int\limits_{n}^{\infty}\frac{\psi^{p'}(t)t^{p'-2}}{\alpha(\psi;t)}dt-(p'-1)\int\limits_{n}^{\infty}\psi^{p'}(t)t^{p'-2}dt.
\end{equation}
Покажемо, що
\begin{equation}\label{for304}
\frac{1}{\alpha(\psi;t)}-\frac{1}{\alpha(g_{p};t)}
=\frac{1}{p}, \ t\geq1, \ \ 1<p<\infty.
\end{equation}
Дійсно,
оскільки
$\psi(t)=g_{p}(t)t^{-\frac{1}{p}}$, то
$$
\frac{1}{\alpha(\psi;t)}=\frac{t|\psi'(t)|}{\psi(t)}=
\frac{-t\big(\frac{-1}{p}t^{-\frac{1}{p}-1}g_{p}(t)+t^{-\frac{1}{p}}g_{p}'(t)\big)}{g_{p}(t)t^{-\frac{1}{p}}}=
\frac{\frac{1}{p}t^{-\frac{1}{p}}g_{p}(t)+t^{-\frac{1}{p}+1}|g_{p}'(t)|}{g_{p}(t)t^{-\frac{1}{p}}}=
$$
$$
=\frac{1}{p}+\frac{t|g_{p}'(t)|}{g_{p}(t)}=\frac{1}{p}+\frac{1}{\alpha(g_{p};t)}.
$$

В силу (\ref{eqq3}) i (\ref{for304})
$$
\psi^{p'}(n)n^{p'-1}=p'\int\limits_{n}^{\infty}\psi^{p'}(t)t^{p'-2}\Big(\frac{1}{p}+\frac{1}{\alpha(g_{p};t)}\Big)dt-(p'-1)\int\limits_{n}^{\infty}\psi^{p'}(t)t^{p'-2}dt=
$$
\begin{equation}\label{eqq5}
=
p'\int\limits_{n}^{\infty}\psi^{p'}(t)t^{p'-2}\frac{1}{\alpha(g_{p};t)}dt.
\end{equation}

З (\ref{k}), (\ref{eqq2}) i (\ref{eqq5})  випливають співвідношення
\begin{equation}\label{for305}
\psi^{p'}(n)n^{p'-1} \leq
\frac{p'}{\underline{\alpha}_{n}(g_{p})}\int\limits_{n}^{\infty}\psi^{p'}(t)t^{p'-2}dt \leq
\frac{p'}{\underline{\alpha}_{n}(g_{p})}
\sum\limits_{k=n}^{\infty}\psi^{p'}(k)k^{p'-2}.
\end{equation}
Нерівність (\ref{lemma_3}) доведено.

Нехай $g_{p}\in\mathfrak{M}_{C}$. Оскільки $\mathfrak{M}_{C}\subset\mathfrak{M}_{0}$, то справедливість другої нерівності в (\ref{lemma_3_1}) випливає з (\ref{lemma_3}).

Врахувавши (\ref{for305}), маємо
$$
\sum\limits_{k=n}^{\infty}\psi^{p'}(k)k^{p'-2}\leq\psi^{p'}(n)n^{p'-2}+\int\limits_{n}^{\infty}\psi^{p'}(t)t^{p'-2}dt\leq
$$
\begin{equation}\label{j2}
 \leq\Big(\frac{p'}{\underline{\alpha}_{n}(g_{p})}\cdot\frac{1}{n}+1\Big)\int\limits_{n}^{\infty}\psi^{p'}(t)t^{p'-2}dt.
\end{equation}

На підставі формул  (\ref{eqq5}) і (\ref{j2}) одержуємо
$$
\psi^{p'}(n)n^{p'-1} \geq
\frac{p'}{\overline{\alpha}_{n}(g_{p})}\int\limits_{n}^{\infty}\psi^{p'}(t)t^{p'-2}dt \geq
\frac{p'}{\overline{\alpha}_{n}(g_{p})}\Big(\frac{p'}{\underline{\alpha}_{n}(g_{p})}\cdot\frac{1}{n}+1\Big)^{-1}
\sum\limits_{k=n}^{\infty}\psi^{p'}(k)k^{p'-2}.
$$

Тим самим співвідношення (\ref{lemma_3_1}), а отже і лему 3, доведено.

З нерівностей (\ref{gg10}), (\ref{eq60}) і (\ref{lemma_3}) отримуємо
 \begin{equation}\label{for307}
  {\cal E}_{n}(C^{\psi}_{\beta,p})_{C}\leq
\frac{1}{\pi}\xi(p')\Big(1+ \frac{p'}{\underline{\alpha}_{n}(g_{p})}\Big)^{\frac{1}{p'}}
\Big(\sum\limits_{k=n}^{\infty}\psi^{p'}(k)k^{p'-2}\Big)^{\frac{1}{p'}}.
 \end{equation}
 Оскільки послідовність ${\underline{\alpha}}_{n}(g_{p})$ монотонно неспадає, то для довільних $n\in\mathbb{N}$
 \begin{equation}\label{eqq7}
  {\cal E}_{n}(C^{\psi}_{\beta,p})_{C}\leq
K_{\psi,p}^{(2)}
\Big(\sum\limits_{k=n}^{\infty}\psi^{p'}(k)k^{p'-2}\Big)^{\frac{1}{p'}}, \ 1<p<\infty,
 \ \frac{1}{p}+\frac{1}{p'}=1,
 \end{equation}
 де
 $$
 K_{\psi,p}^{(2)}=\frac{1}{\pi}\xi(p')\Big(1+ \frac{p'}{\underline{\alpha}_{1}(g_{p})}\Big)^{\frac{1}{p'}}.
 $$

 Знайдемо оцінку знизу величини ${\cal E}_{n}(C^{\psi}_{\beta,p})_{C}$. З цією метою розглянемо  функцію
\begin{equation}\label{eq4}
  f^{*}(t)= f^{*}(\psi;p;n;t)=
\frac{\lambda}{\Big(\sum\limits_{k=n}^{\infty}\psi^{p'}(k)k^{p'-2}\Big)^{\frac{1}{p}}}\sum\limits_{k=n}^{\infty}
\psi^{p'}(k)k^{p'-2}\cos
kt,
\end{equation}
де
\begin{equation}\label{a}
 \lambda=\lambda(\psi;p;n):=\frac{1}{\xi(p)}\Big(\frac{\underline{\alpha}_{n}(g_{p})}{p'+\underline{\alpha}_{n}(g_{p})}\Big)^{\frac{1}{p}},  \ 1<p<\infty,
\ \ \frac{1}{p}+\frac{1}{p'}=1,
\end{equation}
a величини $\xi(p)$  і $\underline{\alpha}_{n}(g_{p})$ означаються формулами (\ref{for100}) i (\ref{k}) відповідно.

З умови $\sum\limits_{k=1}^{\infty}\psi^{p'}(k)k^{p'-2}<\infty$ випливає, що $f^{*}\in C$. Покажемо,  що $\|(f^{*}(t))^{\psi}_{\beta}\|_{p}\leq 1$. Оскільки, згідно з  умовою теореми 1,
$g_{p}\in\mathfrak{M}_{0}$, то функція $\psi^{p'-1}(t)t^{p'-2}=g_{p}^{p'-1}(t)t^{- \ \frac{1}{p'}}$ монотонно спадає до нуля.
В силу означення $(\psi,\beta)$--похідної, майже скрізь виконується рівність
\begin{equation}\label{f}
(f^{*}(t))^{\psi}_{\beta}=\frac{\lambda}{\Bigg(\sum\limits_{k=n}^{\infty}\psi^{p'}(k)k^{p'-2}\Bigg)^{\frac{1}{p}}}\sum\limits_{k=n}^{\infty}
\psi^{p'-1}(k)k^{p'-2}\cos\Big(kt+\frac{\beta\pi}{2}\Big).
\end{equation}
Поклавши в умовах леми 1 $a_{k}=\psi^{p'-1}(k)k^{p'-2}$, $\gamma=\frac{\beta\pi}{2}$,  $s=p$, із (\ref{f}) і (\ref{lemma_3}), одержуємо
$$
\|(f^{*})^{\psi}_{\beta}\|_{p}\leq
$$
$$
\leq \frac{\lambda\xi(p)}{\bigg(\sum\limits_{k=n}^{\infty}\psi^{p'}(k)k^{p'-2}\bigg)^{\frac{1}{p}}}
\bigg(\sum\limits_{k=n}^{\infty}
\big(\psi^{p'-1}(k)k^{p'-2}\big)^{p}k^{p-2}+\psi^{(p'-1)p}(n)n^{(p'-2)p}n^{p-1}\bigg)^{\frac{1}{p}}=
$$
\begin{equation}\label{for2}
=\lambda\xi(p)\bigg(1+ \frac{\psi^{p'}(n)n^{p'-1}}{\sum\limits_{k=n}^{\infty}\psi^{p'}(k)k^{p'-2}}\bigg)^{\frac{1}{p}}\leq
\lambda\xi(p)\Big(1+ \frac{p'}{\underline{\alpha}_{n}(g_{p})}\Big)^{\frac{1}{p}}= 1.
\end{equation}
З (\ref{for2}) випливає   включення
$f^{*}\in C^{\psi}_{\beta,p}$.

Оскільки послідовність ${\underline{\alpha}}_{n}(g_{p})$ монотонно неспадає, то в силу (\ref{a})
$$
\lambda\geq\frac{1}{\xi(p)}\Big(\frac{\underline{\alpha}_{1}(g_{p})}{p'+\underline{\alpha}_{1}(g_{p})}\Big)^{\frac{1}{p}}, \ n\in\mathbb{N},
$$
i для функції $f^{*}$  має
місце оцінка
$$
{\cal E}_{n}(C^{\psi}_{\beta,p})_{C}\geq|f^{*}(0)-S_{n-1}(f^{*},0)|=
$$
$$
=
\frac{\lambda}{\Big(\sum\limits_{k=n}^{\infty}\psi^{p'}(k)k^{p'-2}\Big)^{\frac{1}{p}}}\sum\limits_{k=n}^{\infty}
\psi^{p'}(k)k^{p'-2}=\lambda\bigg(\sum\limits_{k=n}^{\infty}\psi^{p'}(k)k^{p'-2}\bigg)^{\frac{1}{p'}}\geq
$$
\begin{equation}\label{for1010}
\geq
K_{\psi,p}^{(1)}
\bigg(\sum\limits_{k=n}^{\infty}\psi^{p'}(k)k^{p'-2}\bigg)^{\frac{1}{p'}}, \ n\in\mathbb{N},
\end{equation}
де
 $$
 K_{\psi,p}^{(1)}=\frac{1}{\xi(p)}\Big(\frac{\underline{\alpha}_{1}(g_{p})}{p'+\underline{\alpha}_{1}(g_{p})}\Big)^{\frac{1}{p}}.
 $$
 Об'єднуючи (\ref{eqq7}), (\ref{a}) і (\ref{for1010}) отримуємо (\ref{theorem_1}). Теорему 1 доведено.

{\bf Зауваження 1. }{\it В ході доведення теореми 1 за виконання її умов було показано, що для  довільних $n\in\mathbb{N}$ виконується більш точна, порівняно з (\ref{theorem_1}), оцінка:
$$
\frac{1}{\xi(p)}\Big(\frac{\underline{\alpha}_{n}(g_{p})}{p'+\underline{\alpha}_{n}(g_{p})}\Big)^{\frac{1}{p}}  \Big(\sum\limits_{k=n}^{\infty}\psi^{p'}(k)k^{p'-2}\Big)^{\frac{1}{p'}} \leq{\cal E}_{n}(C^{\psi}_{\beta,p})_{C}
\leq \ \ \ \ \ \ \ \ \ \ \ \ \ \ \ \ \ \ \ \ \ \ \ \ \ \
$$
\begin{equation}\label{remark}
\ \ \ \ \ \ \ \ \ \ \ \ \ \ \ \ \ \ \ \ \ \ \ \ \ \ \leq
\frac{1}{\pi}\xi(p')\Big(\frac{p'+\underline{\alpha}_{n}(g_{p})}{\underline{\alpha}_{n}(g_{p})}\Big)^{\frac{1}{p'}}
\Big(\sum\limits_{k=n}^{\infty}\psi^{p'}(k)k^{p'-2}\Big)^{\frac{1}{p'}},
\end{equation}
де $\xi(p),\ 1<p<\infty,$  і $\underline{\alpha}_{n}(g_{p})$--- додатні величини, що означаються за
допомогою формул (\ref{for100}) і (\ref{k}) відповідно.  }

{\bf Теорема 2.} {\it Нехай $\sum\limits_{k=1}^{\infty}\psi^{p'}(k)k^{p'-2}<\infty$,
 $\psi(t)=g_{p}(t)t^{-\frac{1}{p}}$, де $g_{p}\in\mathfrak{M}_{0}$, $1<p<\infty$, $\frac{1}{p}+\frac{1}{p'}=1$ i
 \begin{equation}\label{cond}
  \underline{ \alpha}_{1}(g_{p})=\inf\limits_{t\geq1}\alpha(g_{p};t)>\frac{p'}{2}.
 \end{equation}
  Тоді для довільних $n\in \mathbb{N}$ i $\beta\in\mathbb{R}$ мають місце співвідношення
 \begin{equation}\label{theorem_2}
K^{(3)}_{\psi,p}\Big(\sum\limits_{k=n}^{\infty}\psi^{p'}(k)k^{p'-2}\Big)^{\frac{1}{p'}}\leq{E}_{n}(C^{\psi}_{\beta,p})_{C}
\leq \mathcal{E}_{n}(C^{\psi}_{\beta,p})_{C}\leq
K^{(2)}_{\psi,p}\Big(\sum\limits_{k=n}^{\infty}\psi^{p'}(k)k^{p'-2}\Big)^{\frac{1}{p'}}.
\end{equation}
в яких $K^{(2)}_{\psi,p}$, $K^{(3)}_{\psi,p}$ --- додатні величини, що
  залежать лише  від $\psi$ і $p$.}

{\bf Доведення теореми 2.} В силу теореми 1 якщо
 $\psi(t)=g_{p}(t)t^{-\frac{1}{p}}$, $g_{p}\in\mathfrak{M}_{0}$
 i $\sum\limits_{k=1}^{\infty}\psi^{p'}(k)k^{p'-2}<\infty$, $1<p<\infty$, $\frac{1}{p}+\frac{1}{p'}=1$,
  то існує стала $K^{(2)}_{\psi,p}$ така, що при всіх $n\in \mathbb{N}$ i $\beta\in\mathbb{R}$
\begin{equation}\label{j3}
  {E}_{n}(C^{\psi}_{\beta,p})_{C} \leq \mathcal{E}_{n}(C^{\psi}_{\beta,p})_{C}\leq K^{(2)}_{\psi,p}\Big(\sum\limits_{k=n}^{\infty}\psi^{p'}(k)k^{p'-2}\Big)^{\frac{1}{p'}}.
\end{equation}

 Знайдемо оцінку знизу для величин ${ E}_{n}(C^{\psi}_{\beta,p})_{C}$.
Покладемо
 $$\Phi_{p'}(x):=\int\limits_{x}^{\infty}\psi^{p'}(t)t^{p'-2}dt$$
 і
\begin{equation}\label{for31}
A(n)=A(\psi;p;n):=\big[\Phi_{p'}^{-1}\big(\frac{K^{*}}{n}\Phi_{p'}(n)\big)\big]+1,
\end{equation}
 де $[\alpha]$ --- ціла частина дійсного числа $\alpha$,  $\Phi_{p'}^{-1}$ --- функція обернена до $\Phi_{p'}$,
 a
 \begin{equation}\label{kstar}
K^{*}=K^{*}(\psi;p):=\frac{1}{2}\Big(1-\frac{p'}{2\underline{\alpha}_{1}(g_{p})}\Big), \ \ 1<p<\infty, \ \ \frac{1}{p}+\frac{1}{p'}=1.
\end{equation}

Розглянемо функцію $f^{*}$ вигляду (\ref{eq4}). Як зазначалось раніше (див. (\ref{for2})), $f^{*}\in C^{\psi}_{\beta,p}$, ${1<p<\infty}$.
Покажемо, що має місце оцінка
\begin{equation}\label{j5}
  {E}_{n}(f^{*})_{C}= \inf\limits_{t_{n-1}\in\mathcal{T}_{2n-1}}\|f^{*}(\cdot)-t_{n-1}(\cdot)\|_{C} \geq K^{(3)}_{\psi,p}\Big(\sum\limits_{k=n}^{\infty}\psi^{p'}(k)k^{p'-2}\Big)^{\frac{1}{p'}},
\end{equation}
де $K^{(3)}_{\psi,p}$ --- деяка додатня стала,
 що може залежати лише  від $\psi$ і $p$.

Розглянемо інтеграл
\begin{equation}\label{eq3}
  I=\int\limits_{-\pi}^{\pi}(f^{*}(t)-t_{n-1}(t))(V_{A(n)}(t)-V_{n-1}(t))dt,
\end{equation}
де  $t_{n-1}\in \mathcal{T}_{2n-1}$, а  $V_{m}(t)$ --- ядра Валле Пуссена вигляду
$$
V_{m}(t)=\frac{1}{2}+\frac{1}{m}\sum\limits_{k=m}^{2m-1}\sum\limits_{j=1}^{k}\cos jt, \ m\in
\mathbb{N},
$$
а послідовність $A(n)$ означається формулою (\ref{for31}).

 В силу твердження Д.1.1 з  \cite[с. 391]{Korn}  
\begin{equation}\label{for200}
 I\leq\|f^{*}(t)-t_{n-1}(t)\|_{C}\|V_{A(n)}(t)-V_{n-1}(t)\|_{1}.
\end{equation}
Знайдемо оцінку норми $\|V_{A(n)}(t)-V_{n-1}(t)\|_{1}$.
Відомо, що
\begin{equation}\label{for202}
V_{m}(t)=2F_{2m-1}(t)-F_{m-1}(t),
\end{equation}
 (див., наприклад,  \cite[с. 28]{T}),  де $F_{k}(t)$ --- ядра Фейєра порядку $k$
 $$
 F_{k}(t)=\frac{1}{2}+\frac{1}{k+1}\sum\limits_{\nu=0}^{k}\Big(\sum\limits_{j=1}^{\nu}\cos jt\Big), \ k\in
\mathbb{N}.
 $$
Оскільки (див., наприклад,  \cite[с. 148]{Zigmund1})
\begin{equation}\label{for201}
 \|F_{k}(t)\|_{1}=\pi, \ k\in
\mathbb{N},
\end{equation}
то з (\ref{for202}) i (\ref{for201}) отримуємо
$$
\|V_{m}(t)\|_{1}\leq2\|F_{2m-1}(t)\|_{1}+\|F_{m-1}(t)\|_{1}\leq3\pi, \ \ m\in\mathbb{N}.
$$
Тому
\begin{equation}\label{eq61}
  \|V_{A(n)}(t)-V_{n-1}(t)\|_{1}\leq \|V_{A(n)}(t)\|_{1}+\|V_{n-1}(t)\|_{1}\leq6\pi.
\end{equation}
З (\ref{for200}) і (\ref{eq61}), маємо
\begin{equation}\label{eq9}
 I\leq6\pi\|f^{*}(t)-t_{n-1}(t)\|_{C}.
\end{equation}
Для ядер $V_{m}(t)$
виконується рівність (див., наприклад, формулу
 (1.3.15) роботи  \cite[с. 31]{Stepanets1})
 $$
 V_{m}(t)=\frac{1}{2}+\sum\limits_{k=1}^{m}\cos kt+2\sum\limits_{k=m+1}^{2m-1}\Big(1-\frac{k}{2m}\Big)\cos
kt, \ m\in \mathbb{N},
 $$
а тому
$$
V_{A(n)}(t)-V_{n-1}(t)=
$$
\begin{equation}\label{eq2}
 =\sum\limits_{k=n}^{A(n)}\cos kt+2\sum\limits_{k=A(n)+1}^{2A(n)-1}\Big(1-\frac{k}{2A(n)}\Big)\cos
kt-2\sum\limits_{k=n}^{2n-3}\Big(1-\frac{k}{2n-2}\Big)\cos
kt.
\end{equation}

 Враховуючи очевидні рівності
 $$
 \int\limits_{-\pi}^{\pi}\cos kt\cos mtdt=
{\left\{\begin{array}{cc}
0, \ & k\neq m, \\
\pi, &
k=m, \
  \end{array} \right.}  \ \ k,m\in\mathbb{N}, \
 $$
з формул  (\ref{eq4}), (\ref{eq3}) і (\ref{eq2})  отримаємо, що для довільних $t_{n-1}\in\mathcal{T}_{2n-1}$
$$
I=\int\limits_{-\pi}^{\pi}(f^{*}(t)-t_{n-1})(V_{A(n)}(t)-V_{n-1}(t))dt=
$$
$$
=\int\limits_{-\pi}^{\pi}f^{*}(t)(V_{A(n)}(t)-V_{n-1}(t))dt=
$$
$$
=\frac{\lambda}{\Big(\sum\limits_{k=n}^{\infty}\psi^{p'}(k)k^{p'-2}\Big)^{\frac{1}{p}}}\int\limits_{-\pi}^{\pi}\sum\limits_{k=n}^{\infty}
\psi^{p'}(k)k^{p'-2}\cos
kt\times
$$
$$
\times \bigg( \sum\limits_{k=n}^{A(n)}\cos kt+2\sum\limits_{k=A(n)+1}^{2A(n)-1}\Big(1-\frac{k}{2A(n)}\Big)\cos
kt-2\sum\limits_{k=n}^{2n-3}\Big(1-\frac{k}{2n-2}\Big)\cos
kt \bigg) dt =
$$
$$
=\frac{\pi \lambda}{\Big(\sum\limits_{k=n}^{\infty}\psi^{p'}(k)k^{p'-2}\Big)^{\frac{1}{p}}}\bigg(\sum\limits_{k=n}^{A(n)}
\psi^{p'}(k)k^{p'-2}+2\sum\limits_{k=A(n)+1}^{2A(n)-1}\psi^{p'}(k)k^{p'-2}\Big(1-\frac{k}{2A(n)}\Big)-
$$
$$
- 2\sum\limits_{k=n}^{2n-3}\psi^{p'}(k)k^{p'-2}\Big(1-\frac{k}{2n-2}\Big)\bigg)>
$$
\begin{equation}\label{eq8}
>\frac{\pi \lambda}{\Big(\sum\limits_{k=n}^{\infty}\psi^{p'}(k)k^{p'-2}\Big)^{\frac{1}{p}}}\bigg(\sum\limits_{k=n}^{A(n)}
\psi^{p'}(k)k^{p'-2}
- 2\sum\limits_{k=n}^{2n-3}\psi^{p'}(k)k^{p'-2}\Big(1-\frac{k}{2n-2}\Big)\bigg).
\end{equation}

Знайдемо оцінку зверху для суми $\sum\limits_{k=n}^{2n-3}\psi^{p'}(k)k^{p'-2}\Big(1-\frac{k}{2n-2}\Big)$. Враховуючи нерівність (\ref{lemma_3}) та спадання функції $\psi^{p'}(t)t^{p'-2}$, одержимо
$$
2\sum\limits_{k=n}^{2n-3}\psi^{p'}(k)k^{p'-2}\Big(1-\frac{k}{2n-2}\Big)
\leq\psi^{p'}(n)n^{p'-2}\frac{1}{n-1}\sum\limits_{k=n}^{2n-3}(2n-k-2)
=
$$
\begin{equation}\label{eq5}
 =\psi^{p'}(n)n^{p'-2}\frac{n-2}{2}<\frac{1}{2}\psi^{p'}(n)n^{p'-1}<\frac{p'}{2\underline{\alpha}_{n}(g_{p})}\sum\limits_{k=n}^{\infty}\psi^{p'}(k)k^{p'-2}.
\end{equation}

Беручи до уваги співвідношення (\ref{eq9}), (\ref{eq8}) і (\ref{eq5}), запишемо
$$
\|f^{*}(t)-t_{n-1}(t)\|_{C}\geq \frac{1}{6\pi}I\geq
$$
$$
\geq \frac{1}{6\pi}\frac{\pi \lambda}{\Big(\sum\limits_{k=n}^{\infty}\psi^{p'}(k)k^{p'-2}\Big)^{\frac{1}{p}}}\bigg(\sum\limits_{k=n}^{A(n)}
\psi^{p'}(k)k^{p'-2}
-\frac{p'}{2\underline{\alpha}_{n}(g_{p})}\sum\limits_{k=n}^{\infty}\psi^{p'}(k)k^{p'-2}\bigg)>
$$
\begin{equation}\label{eq6}
>\frac{1}{6}\frac{\lambda}{\Big(\sum\limits_{k=n}^{\infty}\psi^{p'}(k)k^{p'-2}\Big)^{\frac{1}{p}}}
\bigg(\Big(1-\frac{p'}{2\underline{\alpha}_{n}(g_{p})}\Big)\sum\limits_{k=n}^{\infty}\psi^{p'}(k)k^{p'-2}
-\sum\limits_{k=A(n)+1}^{\infty}\psi^{p'}(k)k^{p'-2}\bigg).
\end{equation}

З (\ref{for31}) та монотонного спадання функції $\Phi_{p'}(\cdot)$ випливає, що
$$
\Phi_{p'}(A(n))=\Phi_{p'}(\big[\Phi_{p'}^{-1}\big(\frac{K^{*}}{n}\Phi_{p'}(n)\big)\big]+1)<\frac{K^{*}}{n}\Phi_{p'}(n),
$$
а тому
\begin{equation}\label{eq10}
 \sum\limits_{k=A(n)+1}^{\infty}\psi^{p'}(k)k^{p'-2}\leq\int\limits_{A(n)}^{\infty}\psi^{p'}(t)t^{p'-2}dt=
\Phi_{p'}(A(n))<\frac{K^{*}}{n}\Phi_{p'}(n)\leq\frac{K^{*}}{n}\sum\limits_{k=n}^{\infty}\psi^{p'}(k)k^{p'-2},
\end{equation}
при цьому в силу (\ref{cond}) i (\ref{kstar}) $0<K^{*}<\frac{1}{2}$.

Із (\ref{eq6}) i (\ref{eq10})  з урахуванням монотонного неспадання послідовності $\underline{\alpha}_{n}(g_{p})$, для полінома  $t^{*}_{n-1}(t)$ найкращого рівномірного наближення функції $f^{*}$  отримаємо оцінку
$$
E_{n}(f^{*})_{C}=\|f^{*}(t)-t^{*}_{n-1}(t)\|_{C}\geq
$$
$$
\geq\frac{1}{6}\frac{\lambda}{\Big(\sum\limits_{k=n}^{\infty}\psi^{p'}(k)k^{p'-2}\Big)^{\frac{1}{p}}}
\Big(1-\frac{p'}{2\underline{\alpha}_{n}(g_{p})}-\frac{K^{*}}{n}\Big)
\sum\limits_{k=n}^{\infty}\psi^{p'}(k)k^{p'-2}\geq
$$
\begin{equation}\label{for101}
\geq K_{\psi,p}^{(3)}
\Big(\sum\limits_{k=n}^{\infty}\psi^{p'}(k)k^{p'-2}\Big)^{\frac{1}{p'}},
\end{equation}
де
\begin{equation}\label{k3}
K_{\psi,p}^{(3)}=\frac{1}{12\xi(p)}\Big(\frac{\underline{\alpha}_{1}(g_{p})}{p'+\underline{\alpha}_{1}(g_{p})}\Big)^{\frac{1}{p}}
\Big(1-\frac{p'}{2\underline{\alpha}_{1}(g_{p})}\Big).
\end{equation}

З (\ref{for101})  випливає справедливість оцінки (\ref{j5}).
Об'єднуючи (\ref{j3}) і (\ref{j5}) отримуємо (\ref{theorem_2}).
Теорему 2 доведено.

{\bf Зауваження 2. }{\it З (\ref{remark}) та  ходу доведення теореми 2 випливає, що за виконання умов:
  $\sum\limits_{k=1}^{\infty}\psi^{p'}(k)k^{p'-2}<\infty$,
 $\psi(t)=g_{p}(t)t^{-\frac{1}{p}}$, де $g_{p}\in\mathfrak{M}_{0}$, $1<p<\infty$, $\frac{1}{p}+\frac{1}{p'}=1$, $n\in\mathbb{N}$ і
 $$
 \underline{ \alpha}_{n}(g_{p})=\inf\limits_{t\geq1}\alpha(g_{p};t)>\frac{p'}{2},
 $$
   для  довільного $\beta\in\mathbb{R}$ виконується більш точна, порівняно з (\ref{theorem_2}), оцінка:
$$
\frac{1}{12\xi(p)}\Big(\frac{\underline{\alpha}_{n}(g_{p})}{p'+\underline{\alpha}_{n}(g_{p})}\Big)^{\frac{1}{p}}
\Big(1-\frac{p'}{2\underline{\alpha}_{n}(g_{p})}\Big) \Big(\sum\limits_{k=n}^{\infty}\psi^{p'}(k)k^{p'-2}\Big)^{\frac{1}{p'}} \leq{ E}_{n}(C^{\psi}_{\beta,p})_{C}
 \leq{\cal E}_{n}(C^{\psi}_{\beta,p})_{C}
\leq\ \ \ \ \ \ \ \ \ \ \
$$
\begin{equation}\label{remark2}
\ \ \ \ \ \ \ \ \ \ \ \ \ \ \ \ \ \ \ \ \ \ \ \ \ \ \leq
\frac{1}{\pi}\xi(p')\Big(\frac{p'+\underline{\alpha}_{n}(g_{p})}{\underline{\alpha}_{n}(g_{p})}\Big)^{\frac{1}{p'}}
\Big(\sum\limits_{k=n}^{\infty}\psi^{p'}(k)k^{p'-2}\Big)^{\frac{1}{p'}},
\end{equation}
де $\xi(p),\ 1<p<\infty,$  і $\underline{\alpha}_{n}(g_{p})$--- додатні величини, що означаються формулами (\ref{for100}) і (\ref{k}) відповідно.  }

Надалі, як зазвичай прийнято,  для  послідовностей $A(n)>0$ i $B(n)>0$ під записом  ${A(n)\asymp B(n)}$ будемо розуміти, існування  сталих $K_{1}>0$ і
$K_{2}>0$ таких, що ${K_{1}B(n)\leq A(n)\leq K_{2}B(n)}$,
$n\in\mathbb{N}$.

Із теорем 1 і 2 безпосередньо випливає наступне твердження

{\bf Наслідок 1. }{\it   Нехай $\sum\limits_{k=1}^{\infty}\psi^{p'}(k)k^{p'-2}<\infty$, $\psi(t)=g_{p}(t)t^{-\frac{1}{p}}$, де $g_{p}\in\mathfrak{M}_{0}$, $1<p<\infty$, $\frac{1}{p}+\frac{1}{p'}=1$, $n\in \mathbb{N}$ i $\beta\in\mathbb{R}$.
  Тоді
  \begin{equation}\label{s1}
{\cal E}_{n}(C^{\psi}_{\beta,p})_{C}\asymp\Big(\sum\limits_{k=n}^{\infty}\psi^{p'}(k)k^{p'-2}\Big)^{\frac{1}{p'}}.
\end{equation}
Якщо, крім того, виконується  умова (\ref{cond}), то }
\begin{equation}\label{j7}
{E}_{n}(C^{\psi}_{\beta,p})_{C}\asymp{\cal E}_{n}(C^{\psi}_{\beta,p})_{C}\asymp\Big(\sum\limits_{k=n}^{\infty}\psi^{p'}(k)k^{p'-2}\Big)^{\frac{1}{p'}}.
\end{equation}

Припустимо, що виконуються умови наслідку 1, і крім того, $g_{p}\in\mathfrak{M}_{C}$, тоді як випливає зі співвідношення (\ref{lemma_3_1}),
$$\sum\limits_{k=n}^{\infty}\psi^{p'}(k)k^{p'-2}\asymp\psi^{p'}(n)n^{p'-1}.$$
Отже, має місце твердження.

{\bf Наслідок 2. }{\it  Нехай $\sum\limits_{k=1}^{\infty}\psi^{p'}(k)k^{p'-2}<\infty$, $\psi(t)=g_{p}(t)t^{-\frac{1}{p}}$ де $g_{p}\in\mathfrak{M}_{C}$, $1<p<\infty$, $\frac{1}{p}+\frac{1}{p'}=1$, $n\in \mathbb{N}$ i $\beta\in\mathbb{R}$.
  Тоді
  \begin{equation}\label{s1}
{\cal E}_{n}(C^{\psi}_{\beta,p})_{C}\asymp\psi(n)n^{\frac{1}{p}},
\end{equation}
Якщо, крім того, виконується  умова (\ref{cond}), то }
\begin{equation}\label{s5}
{ E}_{n}(C^{\psi}_{\beta,p})_{C}\asymp{\cal E}_{n}(C^{\psi}_{\beta,p})_{C}\asymp\psi(n)n^{\frac{1}{p}}.
\end{equation}

Справедливість порядкових оцінок (\ref{s1}) i (\ref{s5}) встановлена раніше в \cite{Serdyuk_grabova}.

Зауважимо, що  коли $g_{p}(t)\in\mathfrak{M}_{0}$ і
\begin{equation}\label{lim}
\lim\limits_{t\rightarrow\infty}\alpha(g_{p};t)=\infty,
\end{equation}
 то порядкові оцінки (\ref{s1}) i (\ref{s5}) місця не мають, оскільки в цьому випадку виконується оцінка
$$
\psi(n)n^{\frac{1}{p}}=o\bigg(\Big(\sum\limits_{k=n}^{\infty}\psi^{p'}(k)k^{p'-2}\Big)^{\frac{1}{p'}}\bigg), \ n\rightarrow\infty,
$$
яка є наслідком нерівності (\ref{lemma_3}) і того, що $\underline{\alpha}_{n}(g_{p}(t))\rightarrow\infty$ при $n\rightarrow\infty$.

 Прикладом функцій $\psi(t)$, які задовольняють умови наслідку 1 і для яких виконується умова (\ref{lim})  є функції
 \begin{equation}\label{func}
\psi(t)={t^{-\frac{1}{p}}\ln^{-\gamma}(t+K)}, \ {\gamma>\frac{1}{p'}},  \ K\geq e^{\frac{\gamma p'}{2}}.
\end{equation}
  Для них
 $$\sum\limits_{k=1}^{\infty}\psi^{p'}(k)k^{p'-2}=\sum\limits_{k=1}^{\infty}\frac{1}{k\ln^{\gamma p'}(k+K)}<\infty$$ i
$$
g_{p}(t)=\ln^{-\gamma}(t+K), \ \
g_{p}'(t)=-\gamma\ln^{-\gamma-1}(t+K)\frac{1}{t+K},
$$
$$
\alpha(g_{p};t)=\frac{\ln(t+K)}{\gamma}\frac{t+K}{t}
>\frac{\ln(t+e^{\frac{\gamma p'}{2}})}{\gamma},
$$
а тому $\alpha(g_{p};t)\rightarrow\infty$ при $t\rightarrow\infty$ i
\begin{equation}\label{for403}
\underline{\alpha}_{1}(g_{p})>\frac{ p'}{2}.
\end{equation}

Наведемо порядкові оцінки величин ${ E}_{n}(C^{\psi}_{\beta,p})_{C}$ i ${\cal E}_{n}(C^{\psi}_{\beta,p})_{C}$  у випадку коли $\psi(t)$ мають вигляд (\ref{func}).

{\bf Наслідок 3. }{\it Нехай $\psi(t)={t^{-\frac{1}{p}}\ln^{-\gamma}(t+K)}$, ${\gamma>\frac{1}{p'}}$, $K\geq e^{\frac{\gamma p'}{2}}$, $1<p<\infty$, $\frac{1}{p}+\frac{1}{p'}=1$,  $\beta\in\mathbb{R}$ i $n\in\mathbb{N}$. Тоді
$$
{ E}_{n}(C^{\psi}_{\beta,p})_{C}\asymp {\cal E}_{n}(C^{\psi}_{\beta,p})_{C}\asymp\psi(n)n^{\frac{1}{p}}\ln^{\frac{1}{p'}} n, \ \  \ n\geq2.
$$

}
{\bf Доведення наслідку 3. }
З (\ref{j2}) i (\ref{for403}), одержуємо
\begin{equation}\label{for110}
\int\limits_{n}^{\infty}\psi^{p'}(t)t^{p'-2}dt
\leq\sum\limits_{k=n}^{\infty}\psi^{p'}(k)k^{p'-2}
\leq3\int\limits_{n}^{\infty}\psi^{p'}(t)t^{p'-2}dt.
\end{equation}
Згідно з  (\ref{j7}) і  (\ref{for110}) для зазаначених $\psi$
$$
{ E}_{n}(C^{\psi}_{\beta,p})_{C}\asymp {\cal E}_{n}(C^{\psi}_{\beta,p})_{C}
 \asymp\Big(\sum\limits_{k=n}^{\infty}\psi^{p'}(k)k^{p'-2}\Big)^{\frac{1}{p'}}
 \asymp\Big(\int\limits_{n}^{\infty}\psi^{p'}(t)t^{p'-2}dt\Big)^{\frac{1}{p'}}=
$$
$$
=\Big(\int\limits_{n}^{\infty}\frac{dt}{t\ln^{\gamma p'}(t+K)}\Big)^{\frac{1}{p'}}
\asymp
 \ln^{\frac{1}{p'}-\gamma} n=\psi(n)n^{\frac{1}{p}}\ln^{\frac{1}{p'}} n\frac{\ln^{-\gamma}(n)}{\ln^{-\gamma}(n+K)}
 \asymp \psi(n)n^{\frac{1}{p}}\ln^{\frac{1}{p'}} n,
 \ n\geq2.
$$
Наслідок 3 доведено.

E-mail: \href{mailto:serdyuk@imath.kiev.ua}{serdyuk@imath.kiev.ua},
\href{mailto:tania_stepaniuk@ukr.net}{tania$_{-}$stepaniuk@ukr.net}


\newpage
\begin{thebibliography}{10}


\bibitem{Stepanets1}{\sc Степанец А.И.}
 \emph{Методы теории приближений}: В 2 ч. // Праці Інституту математики НАН України --- Киев: Ин-т
математики НАН Украины, 2002. --- {\bf 40}. --- Ч.І. --- 427 с.

\bibitem{Zigmund2}{\sc Зигмунд А.}
 \emph{Тригонометрические ряды.}  В 2 т.~--- М.: Мир,
1965.~--- Т. ІІ.~ --- 538 с.


\bibitem{Step monog 1987}{\sc Степанец А.И.}
 \emph{Классификация и приближение периодических функций.} --- Киев: Наук. думка~--- 1987.~--- 268 c.
 
 


\bibitem{Stepanets_Serdyuk_Shydlich}{\sc Степанец А.И., Сердюк А.С., Шидлич А.Л.}
Классификация бесконечно дифференцируемых функций
// Укр. мат. журн. --- 2008.
--- {\bf 60}, №12. --- С. 1686--1708.



\bibitem{T}{\sc Temlyakov V.N.}
 \emph{Approximation of Periodic Function}: Nova Science Publi--\ chers, Inc. --- 1993. --- 419p.
 

\bibitem{Kol}{\sc Kolmogoroff A.}
Zur Gr\"{o}ssennordnung des Restgliedes Fourierschen Reihen
differenzierbarer Funktionen // Ann. Math.(2),
--- 1935. --- {\bf 36}, №2. --- P.~ 521--526.



\bibitem{Pin}{\sc Пинкевич В.Т.}
О порядке остаточного члена ряда Фурье функций, дифференцируемых в
смысле Вейля // Изв. АН СССР. Сер. мат. --- 1940.
--- {\bf 4}, №6. --- С. 521--528.


\bibitem{Nik2}{\sc Никольский С.М.}
Приближение  функций тригонометрическими полиномами в среднем
// Изв. АН
СССР. Cер. мат.
--- 1946. --- {\bf 10}, №3. --- С. 207--256.


\bibitem{Fav}{\sc Favard J.}
Sur l'approximation des fonctions p\'{e}riodiques par des polynomes
trigonom\'{e}triques
// C.R. Acad. Sci. --- 1936. --- {\bf 203}. --- P. 1122--1124.

\bibitem{Dz}{\sc Дзядык В.К.}
О наилучшем приближении на классе периодических функций, имеющих
ограниченную $s$--ю производную $(0<s<1)$ // Изв. АН СССР, Сер. мат.
--- 1953. --- {\bf 17}. --- С. 135--162.

\bibitem{Dz1}{\sc Дзядык В.К.}
О наилучшем приближении на классах периодических функций,
определяемых интегралами от линейной комбинации абсолютно монотонных
ядер  // Мат. заметки. --- 1974. --- {\bf 16}, №5.
--- С.~ 691--701.


\bibitem{ST}{\sc Стечкин С.Б.}
О наилучшем приближении некоторых классов периодических функций
тригонометрическими полиномами // Изв. АН СССР, Cер. мат. --- 1956.
--- {\bf 20}, --- С. 643--648.

\bibitem{SU}{\sc  Сунь Юн--шен.}
О наилучшем приближении периодических дифференцируемых функций
тригонометрическими полиномами // Изв. АН СССР. Cер. мат.
--- 1959. --- {\bf 23}, №1. --- С. 67--92.

\bibitem{Serdyuk}{\sc  Сердюк А.С., Соколенко І.В.}
Рівномірні наближення класів $(\psi,
\overline{\beta})$--диференційовних функцій лінійними методами
// Теорія наближення функцій та суміжні питання: Зб. праць Ін--ту математики НАН України.2011. --- {\bf 8}, №1.
 --- С. 181--189.
 
 \bibitem{Serdyuk_grabova}{\sc  Сердюк А.С., Грабова У.З.}
 Порядкові оцінки найкращих
наближень і наближень сумами Фур'є  класів $(\psi,\beta)$ --
диференційовних функцій
// Укр. мат. журн. --- 2013.
--- {\bf 65}, №9. --- С. 1186 -- 1197.

\bibitem{Rom}{\sc Романюк В.С.}
Дополнения к оценкам приближения суммами Фурье классов бесконечно
дифференцируемых функций // Екстремальні задачі теорії функцій та
суміжні питання: Праці Ін-ту математики НАН України.
--- 2003. --- {\bf 46},
 --- С. 131--135.
 

\bibitem{S_S}{\sc  Сердюк А.С., Степанюк Т.А.}
Порядкові оцінки найкращих
наближень і наближень сумами Фур'є класів нескінченно
диференційовних функцій// Теорія наближення функцій та суміжні питання: Зб. праць Ін--ту математики НАН України.
2013. --- {\bf 10}, №1.
 --- С. 255--282.


\bibitem{Stepanets2}{\sc Степанец А.И.}
 \emph{Методы теории приближений.} --- Киев: Ин-т математики НАН Украины, 2002. --- {\bf 40}.
--- Ч.II. --- 468 с.


\bibitem{Gradshteyn}{\sc Градштейн И.С., Рыжик И.М.}
 \emph{Таблицы интегралов, сумм, рядов и произведений.}  --- М.: Физматиз, 1962. --- 1100 с.

\bibitem{Zigmund1}{\sc Зигмунд А.}
 \emph{Тригонометрические ряды.}  В 2 т.~--- М.: Мир,
1965.~--- Т. І.~ --- 615 с.

\bibitem{Korn}{\sc Корнейчук Н.П.}
 \emph{Точные константы в теории приближения.}  --- М.:
Наука, 1987. ---  424~с.






\end{thebibliography}
\end{document}